\documentclass[leqno,12pt]{article}
\usepackage{amsmath,amssymb,rotating,a4wide,color}
\usepackage{wasysym}
\usepackage{url}
\usepackage{multirow}

\usepackage[all,cmtip]{xy}

\usepackage{verbatim}
\usepackage{hyperref}
\hypersetup{pdfauthor={Nicolas M\"uller, Richard Pink}, pdftitle={Hyperelliptic Curves with Many Automorphisms}}

\newdir^{ (}{{}*!/-3pt/\dir^{(}}
\newdir_{ (}{{}*!/-3pt/\dir_{(}}

\entrymodifiers={+!!<0pt,\fontdimen22\textfont2>}

\def\otau{{\overline{\tau}}}
\def\oA{{\bar{A}}}

\def\oG{{\bar{G}}}
\def\oH{{\bar{H}}}
\def\oK{{\bar{K}}}
\def\oX{{\bar{X}}}
\def\Branch{\mathop{\rm Br}\nolimits}
\def\Jac{\mathop{\rm Jac}\nolimits}
\def\Fix{\mathop{\rm Fix}\nolimits}
\def\triv{{\rm triv}}

\def\bigtimes{\mathop{\raise-2pt\hbox{\huge$\times$}}}

\newbox\circbulletbox
\setbox\circbulletbox=\hbox{\,{\large$\circledcirc$}\kern-8.7pt\raise .6pt\hbox{$\bullet$}\,}

\let\le\leqslant
\let\ge\geqslant

\def\mycirc{{\kern1pt\circ\kern2pt}}

\def\Aut{\mathop{\rm Aut}\nolimits}
\def\Frob{\mathop{\rm Frob}\nolimits}

\def\Gal{\mathop{\rm Gal}\nolimits}
\def\End{\mathop{\rm End}\nolimits}

\def\deg{\mathop{\rm deg}\nolimits}

\def\trace{\mathop{\rm tr}\nolimits}

\def\GL{\mathop{\rm GL}\nolimits}
\def\SL{\mathop{\rm SL}\nolimits}

\def\PGL{\mathop{\rm PGL}\nolimits}

\def\Sym{\mathop{\rm Sym}\nolimits}

\def\der{{\rm der}}

\let\phi\varphi
\let\theta\vartheta
\let\epsilon\varepsilon
\let\setminus\smallsetminus

\newtheorem{Thm}{Theorem}[section]
\newtheorem{Prop}[Thm]{Proposition}

\newtheorem{Cor}[Thm]{Corollary}

\newtheorem{Rem}[Thm]{Remark}

\numberwithin{Thm}{section}
\def\UseTheoremCounterForNextEquation{\setcounter{equation}{\value{Thm}}\addtocounter{Thm}{1}}

\def\qed{{\hskip0pt\unskip\unskip\nobreak\hfil\penalty50
          \hskip1em\hbox{}\nobreak\hfil
           {$\square$}
          \parfillskip=0pt\finalhyphendemerits=0
          \par}\medskip}

\newenvironment{Proof}
               {\noindent{\bf Proof.}\ }
               {\qed}

               {\noindent{\bf Proof of #1.}\ }
               {\qed}

\newcommand{\BA}{{\mathbb{A}}}

\newcommand{\BC}{{\mathbb{C}}}

\newcommand{\BP}{{\mathbb{P}}}
\newcommand{\BQ}{{\mathbb{Q}}}

\newcommand{\BZ}{{\mathbb{Z}}}

\newcommand{\Fp}{{\mathfrak{p}}}

\newcommand{\CA}{{\cal A}}

\newbox\mybox
\def\arrover#1{\mathrel{
       \setbox\mybox=\hbox spread 1.4em
              {\hfil$\scriptstyle#1$\hfil}
       \vbox{\offinterlineskip\copy\mybox
             \hbox to\wd\mybox{\rightarrowfill}}}}

\def\larrover#1{\mathrel{
       \setbox\mybox=\hbox spread 1.4em
              {\hfil$\scriptstyle#1\vphantom{g}$\hfil}
       \vbox{\offinterlineskip\copy\mybox
             \hbox to\wd\mybox{\leftarrowfill}}}}

\def\ontoover#1{\mathrel{
       \setbox\mybox=\hbox spread 1.4em
              {\hfil$\scriptstyle#1\vphantom{g}$\hfil}
       \vbox{\offinterlineskip\copy\mybox
             \hbox to\wd\mybox{\rightarrowfill\hskip-2.8mm
                               $\rightarrow$}}}}
\def\leftontoover#1{\mathrel{
       \setbox\mybox=\hbox spread 1.4em
              {\hfil$\scriptstyle#1\vphantom{g}$\hfil}
       \vbox{\offinterlineskip\copy\mybox
             \hbox to\wd\mybox{$\leftarrow$\hskip-2.8mm
                               \leftarrowfill}}}}

\let\into\hookrightarrow
\let\onto\twoheadrightarrow

\begin{document}

\title{Hyperelliptic Curves with Many Automorphisms}


\author{
\begin{minipage}{.3\hsize}
Nicolas M\"uller\\[12pt]
\small Department of Mathematics \\
ETH Z\"urich\\
8092 Z\"urich\\
Switzerland \\
nicolas.mueller@math.ethz.ch\\[9pt]
\end{minipage}
\qquad
\begin{minipage}{.3\hsize}
Richard Pink\\[12pt]
\small Department of Mathematics \\
ETH Z\"urich\\
8092 Z\"urich\\
Switzerland \\
pink@math.ethz.ch\\[9pt]
\end{minipage}
}
\date{November 17, 2017}
\maketitle

\centerline{To Frans Oort}
\bigskip\bigskip

\begin{abstract}
We determine all complex hyperelliptic curves with many automorphisms and decide which of their jacobians have complex multiplication.
\end{abstract}

{\renewcommand{\thefootnote}{}
\footnotetext{MSC classification: 14H45 (14H37, 14K22)}
}


\newpage


\section{Introduction}
\label{Intro}

Let $X$ be a smooth connected projective algebraic curve of genus $g\ge2$ over the field of complex numbers. Following Rauch \cite{Rauch1970} and Wolfart \cite{Wolfart1997} we say that $X$ has many automorphisms if it cannot be deformed non-trivially together with its automorphism group.
Given his life-long interest in special points on moduli spaces, Frans Oort \cite[Question\;5.18.(1)]{Oort2013} asked whether the point in the moduli space of curves associated to a curve $X$ with many automorphisms is special, i.e., whether the jacobian of $X$ has complex multiplication.

Here we say that an abelian variety $A$ has \emph{complex multiplication} over a field $K$ if $\End_K^\circ(A)$ contains a commutative, semisimple $\BQ$-subalgebra of dimension $2\dim A$. (This property is called ``sufficiently many complex multiplications'' in Chai, Conrad and Oort \cite[Def.\;1.3.1.2]{ChaiConradOort2014}.)

Wolfart \cite{Wolfart2000} observed that the jacobian of a curve with many automorphisms does not generally have complex multiplication and answered Oort's question for all $g\le 4$. In the present paper we answer Oort's question for all hyperelliptic curves with many automorphisms.

\medskip
For this we first determine all isomorphism classes of such curves. For any hyperelliptic curve $X$ over $\BC$ the automorphism group $G$ is an extension of degree $2$ of a finite subgroup $\oG$ of $\PGL_2(\BC)$, called the reduced automorphism group of~$X$. The classification proceeds by going through all possibilities for~$\oG$, using a criterion of Wolfart \cite{Wolfart1997}. The result is that the isomorphism classes of hyperelliptic curves with many automorphisms fall into three infinite families with $\oG$ cyclic or dihedral  and 15 further curves with $\oG\cong A_4$, $S_4$, $A_5$.

All this is essentially known: For the infinite families see for instance Wolfart \cite[\S6.1]{Wolfart2000}; for the others see 
Shaska \cite{Shaska2006b}; and we do use the explicit equations from \cite[Table\;2]{Shaska2006b}. But as we restrict ourselves strictly to hyperelliptic curves with many automorphisms, we can present the classification more succinctly.

\medskip
A list of all hyperelliptic curves with many automorphisms up to isomorphism is collated in Table~\ref{TableHyperellipticManyAuto}. The equations are given in terms of certain separable polynomials from Table~\ref{TablePolybomials}. For the sake of completeness Table~\ref{TableHyperellipticManyAuto} also contains a description of $G=\Aut_\BC(X)$ in all cases, taken from and in the notation of Shaska \cite[Table\;1]{Shaska2006b}.

\begin{table}[h]
\[\begin{array}{|c||c|c|c|c|c|}
\hline
{\large\strut} X      & \oG		&\text{Genus}	&\text{Affine equation} & G					& \Jac(X)\text{ has}\\
\hline\hline
{\large\strut} X_1    & C_{2g+1}	& g\ge 2 		& y^2=x^{2g+1}-1		& C_{4g+2}			& \text{CM}\\\hline
{\large\strut} X_2    & D_{2g+2}	& g\ge 2 		& y^2=x^{2g+2}-1 		& V_{2g+2}			& \text{CM}\\\hline
{\large\strut} X_3    & D_{2g}		& g\ge 3 		& y^2=x^{2g+1}-x 		& U_{2g}			& \text{CM}\\\hline
{\large\strut} X_4    & A_4		& 4 			& y^2=t_4p_4 			& \SL_2(3)			& \text{CM}\\\hline
{\large\strut} X_5    & S_4		& 2				& y^2=t_4 				& \GL_2(3)			& \text{CM}\\\hline
{\large\strut} X_6    & S_4		& 3				& y^2=s_4 				& C_2\times S_4		& \text{no CM}\\\hline
{\large\strut} X_7    & S_4		& 5				& y^2=r_4 				& W_2				& \text{CM}\\\hline
{\large\strut} X_8    & S_4		& 6				& y^2=s_4t_4 			& \GL_2(3)			& \text{no CM}\\\hline
{\large\strut} X_9    & S_4		& 8				& y^2=r_4t_4 			& W_3 				& \text{CM}\\\hline
{\large\strut} X_{10} & S_4		& 9				& y^2=r_4s_4 			& W_2				& \text{no CM}\\\hline
{\large\strut} X_{11} & S_4		& 12			& y^2=r_4s_4t_4			& W_3				& \text{no CM}\\\hline
{\large\strut} X_{12} & A_5		& 5				& y^2=s_5 				& C_2 \times A_5	& \text{no CM}\\\hline
{\large\strut} X_{13} & A_5		& 9				& y^2=r_5 				& C_2 \times A_5 	& \text{no CM}\\\hline
{\large\strut} X_{14} & A_5		& 14			& y^2=t_5				& \SL_2(5)			& \text{CM}\\\hline
{\large\strut} X_{15} & A_5		& 15			& y^2=r_5s_5			& C_2 \times A_5 	& \text{no CM}\\\hline
{\large\strut} X_{16} & A_5		& 20			& y^2=s_5t_5 			& \SL_2(5)			& \text{no CM}\\\hline
{\large\strut} X_{17} & A_5		& 24			& y^2=r_5t_5			& \SL_2(5)			& \text{no CM}\\\hline
{\large\strut} X_{18} & A_5		& 30			& y^2=r_5s_5t_5			& \SL_2(5)			& \text{no CM}\\\hline
\end{array}\]
\caption{All hyperelliptic curves with many automorphisms}
\label{TableHyperellipticManyAuto}
\end{table}

\begin{table}[h]
$$\begin{array}{|l||c|}
\hline
{\large\strut} t_4 & x(x^4-1) \\
\hline
{\large\strut} p_4 & x^4+2i\sqrt{3}x^2+1 \\
\hline
{\large\strut} q_4 & x^4-2i\sqrt{3}x^2+1 \\
\hline
{\large\strut} r_4 & x^{12} - 33 x^{8} - 33 x^{4} + 1 \\
\hline
{\large\strut} s_4=p_4q_4 & x^{8} + 14 x^{4} + 1 \\
\hline
{\large\strut} r_5 & x^{20} - 228  x^{15} + 494  x^{10} + 228  x^{5} + 1 \\
\hline
{\large\strut} s_5 & x(x^{10} + 11  x^{5} - 1) \\
\hline
{\large\strut} t_5 & x^{30} + 522  x^{25} - 10005  x^{20} - 10005  x^{10} - 522  x^{5} + 1 \\
\hline
\end{array}$$
\caption{Certain separable polynomials over $\BC$}
\label{TablePolybomials}
\end{table}

\medskip
For every curve $X$ in the three infinite families the jacobian has complex multiplication, because $X$ is a quotient of a Fermat curve: see Wolfart \cite[\S6.1]{Wolfart2000}. 

\medskip
For 5 of the other curves the jacobian also has complex multiplication. We establish this by verifying a representation theoretic sufficient condition given by Streit \cite{Streit2001}, which essentially shows that $\Jac(X)$ cannot be deformed non-trivially as a polarized abelian variety together with the action of~$G$. 

\medskip
For each of the 10 remaining curves $X$ it turns out that $\Jac(X)$ does not have complex multiplication. To prove this it suffices to exhibit an abelian subvariety of $\Jac(X)$ without complex multiplication. A natural candidate for this is the jacobian of the quotient of $X$ by a subgroup $H<\Aut_\BC(X)$, whose genus is positive but small. In 5 of the cases we found a quotient $H\backslash X$ of genus~$1$ 
and were done when its $j$-invariant was not an algebraic integer. 

\medskip
In the last 5 cases we only found quotients of genus $2$, $4$, or $6$ (except for a quotient of genus $1$ of $X_{10}$ which does have complex multiplication).
In these cases we first tried to find a place where $\Jac(X)$ has partially multiplicative reduction, using the theory of Bosch \cite[Th.\,4.2]{Bosch1980} that describes the reduction of a hyperelliptic curve at a place of odd residue characteristic. For more details about this see Section 10 of the master thesis of the first author \cite{Mueller2017}, on which much of the present paper is based.
But in all these cases we only found good reduction, and an analogous description of the reduction of a hyperelliptic curve at a place of residue characteristic $2$ is not available.

Instead we formulate and implement a simple criterion for complex multiplication that is based solely on the characteristic polynomials of Frobenius. It relies on the Tate conjecture for endomorphisms of abelian varieties and the fact that a non-trivial semisimple algebraic group over $\BQ_\ell$ always possesses non-isomorphic maximal tori. Thus if $\Jac(X)$ does not have complex multiplication, the characteristic polynomials of Frobenius cannot all split over the same number field. For precise statements see Theorem \ref{CMCritThm} and its corollaries. 
In each of the last 5 cases, we verified this criterion by a quick computation that boiled down to using the characteristic polynomials of Frobenius for at most three primes.

\medskip
All the calculations are performed with computer algebra systems. To find equations for the quotient curves $H\backslash X$ and to verify the criterion about characteristic polynomials of Frobenius we employ Sage \cite{sage}. To verify Streit's representation theoretic criterion we use GAP \cite{GAP4}. The respective worksheets can be downloaded from \cite{MuellerPinkWorksheets} both as text files and as pdf files with output.


\section{A criterion of Wolfart}
\label{Wolf}

Throughout the following we consider a smooth connected projective algebraic curve $X$ of genus $g\ge2$ over~$\BC$ and abbreviate $G:=\Aut_\BC(X)$. Following Rauch \cite{Rauch1970} and Wolfart \cite{Wolfart1997}, \cite{Wolfart2000} we say that $X$ has \emph{many automorphisms} 
if the corresponding point $p$ on the moduli space $M_g$ of compact Riemann surfaces of genus $g$ has (in the complex topology) a neighbourhood $U\subset M_g$ such that the Riemann surface corresponding to any point of $U\setminus\{p\}$ has an automorphism group strictly smaller than~$G$. In other words, the number of automorphisms strictly decreases under proper deformations of~$X$.
 

We will use the following criterion of Wolfart: 

\begin{Thm}\label{WolfartThm}
The following are equivalent:
\begin{enumerate}
\item[(a)] The curve $X$ has many automorphisms.
\item[(b)] There exists a subgroup $H<G$, such that $H\backslash X$ has genus $0$ and the projection morphism $X\onto H\backslash X$ has at most three branch points in $H\backslash X$.
\item[(c)] The quotient $G\backslash X$ has genus $0$ and the projection morphism $X\onto G\backslash X$ has at most three branch points in $G\backslash X$.
\end{enumerate}
\end{Thm}

\begin{Proof}
The condition on $H$ in (b) means that $X\onto H\backslash X\cong\BP^1_\BC$ is a Belyi function. Thus (b) itself is equivalent to saying that there exists a Belyi function $X\to\BP^1_\BC$ defining a normal covering. The equivalence of (a) and (b) is therefore precisely the content of Wolfart \cite[Thm.\;6]{Wolfart1997}. But the proof of \cite[Lemma 8]{Wolfart1997} actually shows that (a) implies (c). Since  (c) trivially implies (b), all three statements are equivalent.
\end{Proof}


\section{Hyperelliptic curves}
\label{Hyp}

By definition $X$ is hyperelliptic if and only if it there exists a morphism $\pi\colon X\to\BP^1_\BC$ of degree~$2$. In that case $\pi$ is a Galois covering and $X$ is determined up to isomorphism by the set $\Branch(\pi)\subset\BP^1_\BC$ of $2g+2$ branch points of~$\pi$. Conversely, for any set of $2g+2$ closed points in $\BP^1_\BC$ there is a hyperelliptic curve with precisely these branch points. 
Moreover, the covering involution $\sigma$ of $\pi$ lies in the center of~$G:=\Aut_\BC(X)$, and the factor group $\oG := G/\langle\sigma\rangle$, called the \emph{reduced automorphism group of~$X$}, embeds into $\Aut_\BC(\BP^1_\BC)\cong\PGL_2(\BC)$. Since $X$ is determined by $\Branch(\pi)$, it turns out (see for instance \cite[\S2]{SevillaShaska2007}) that 
\UseTheoremCounterForNextEquation
\begin{equation}\label{oGFormula}
\oG\ =\ \bigl\{ f\in\PGL_2(\BC) \bigm| f(\Branch(\pi))=\Branch(\pi) \bigr\}.
\end{equation}
In particular $\Branch(\pi)$ is a union of $\oG$-orbits.

\medskip
Let $\Fix(\oG)$ denote the set of closed points in $\BP^1_\BC$ on which $\oG$ does not act freely. Then $\oG\backslash\Fix(\oG) \subset \oG\backslash\BP^1_\BC$ is precisely the set of branch points of the projection morphism $\BP^1_\BC\onto\oG\backslash\BP^1_\BC$. 
Thus the set of branch points of the projection morphism $X\onto G\backslash X\cong \oG\backslash\BP^1_\BC$ is precisely $\oG\backslash(\Branch(\pi)\cup\Fix(\oG))$. Since the quotient $G\backslash X \cong \oG\backslash\BP^1_\BC$ automatically has genus~$0$, Theorem \ref{WolfartThm} shows that $X$ has many automorphisms if and only if the cardinality of $\oG\backslash(\Branch(\pi)\cup\Fix(\oG))$ is at most~$3$.
As we have assumed that $X$ has genus $g\ge2$, and every covering of $\BP^1_\BC$ with fewer than $3$ branch points has genus~$0$, the cardinality must actually be equal to~$3$. 

\medskip
The well-known classification following Klein \cite{Klein1888} and Blichfeldt \cite[\S\S52--55]{Blichfeldt1917} states that every finite subgroup of $\PGL_2(\BC)$ is isomorphic to precisely one of the cyclic group $C_n$ of order $n\ge 1$, the dihedral group $D_n$ of order $2n$ for $n\ge2$, or of $A_4$, $S_4$, $A_5$, and that each isomorphism class of such groups corresponds to precisely one conjugacy class of subgroups of $\PGL_2(\BC)$. The classification also tells us the branch points of $\BP^1_\BC\onto\oG\backslash\BP^1_\BC$. In particular, by Brandt and Stichtenoth \cite[\S2]{BrandtStichtenoth1986} we have
\UseTheoremCounterForNextEquation
\begin{equation}\label{BranchFormula}
\bigl|\oG\backslash\!\Fix(\oG)\bigr|\ =\ \left\{\begin{array}{ll}
0 & \hbox{if $\oG=1$,}\\[3pt]
2 & \hbox{if $\oG\cong C_n$ for $n>1$,}\\[3pt]
3 & \hbox{if $\oG\cong D_n$ for $n>1$ or $A_4$, $S_4$, $A_5$.}
\end{array}\right.
\end{equation}
Combining this with the above criterion we deduce that $X$ has many automorphisms if and only if 
\UseTheoremCounterForNextEquation
\begin{equation}\label{ManyAutCondition}
\left\{\begin{array}{ll}
|\Branch(\pi)|=3 & \hbox{if $\oG=1$,}\\[3pt]
|\oG\backslash(\Branch(\pi)\setminus\Fix(\oG))|=1 & \hbox{if $\oG\cong C_n$ for $n>1$,}\\[3pt]
\Branch(\pi)\subset\Fix(\oG) & \hbox{if $\oG\cong D_n$ for $n>1$ or $A_4$, $S_4$, $A_5$.}
\end{array}\right.
\end{equation}
Since $|\Branch(\pi)|=2g+2$ must be even, the first of these cases is in fact impossible. The second case amounts to saying that $\Branch(\pi)$ contains precisely one free $\oG$-orbit, and the third to saying that $\Branch(\pi)$ contains no free $\oG$-orbit.



\medskip
We can now compile an explicit list of all hyperelliptic curves with many automorphisms. For each conjugacy class of finite subgroups $\oG<\PGL_2(\BC)$ we choose the representative in the coordinates from Shaska \cite{Shaska2006b}. By the above we must have $\oG\not=1$; so we now assume that $n>1$ in both the cases $\oG\cong C_n,D_n$.

For each $\oG$ we first write down the decomposition of $\Fix(\oG)$ into $\oG$-orbits. For this let $V(p)$ denote the set of zeros of a polynomial $p\in\BC[x]$, viewed as a subset of $\BP^1_\BC$ via the usual embedding $\BA^1_\BC\subset\BP^1_\BC$. With the separable polynomials from Table~\ref{TablePolybomials} above, we can extract from \cite[\S4]{Shaska2006b} the following description of $\oG$-orbits:
$$\begin{array}{|c||c|c|}
\hline
{\large\strut} \oG & \hbox{$\oG$-orbits in $\Fix(\oG)$}  & \hbox{respective orbit sizes} \\
\hline\hline
{\large\strut} C_n & \{\infty\},\ \{0\} & 1,\ 1 \\
\hline
{\large\strut} D_n & \{0,\infty\},\ V(x^n-1),\ V(x^n+1) & 2,\ n,\ n  \\
\hline
{\large\strut} A_4 & V(t_4)\cup\{\infty\},\ V(p_4),\ V(q_4) & 6,\ 4,\ 4 \\
\hline
{\large\strut} S_4 & V(t_4)\cup\{\infty\},\ V(r_4),\ V(s_4) & 6,\ 12,\ 8 \\
\hline
{\large\strut} A_5 & V(s_5)\cup\{\infty\},\ V(r_5),\ V(t_5) & 12,\ 20,\ 30 \\
\hline
\end{array}$$



\medskip
Next we list all the possibilities for the branch locus $\Branch(\pi)\subset\BP^1_\BC$. This must be a subset of even cardinality $2g+2\ge6$ which is a union of $\oG$-orbits and satisfies the condition in (\ref{ManyAutCondition}), but is subject to no other requirements.

In the case $\oG\not\cong C_n$ it must be a union of $\oG$-orbits in $\Fix(\oG)$. In the cases $\oG\cong S_4,A_5$ any non-empty union of $\oG$-orbits in $\Fix(\oG)$ is okay, yielding $7$ possibilities for $X$ each. As each of $S_4$, $A_5$ is a maximal finite subgroup of $\PGL_2(\BC)$, the reduced automorphism group of $X$ is then really~$\oG$. Since the orbit structures are different in all these cases, the resulting curves $X$ are pairwise non-isomorphic.

For the case $\oG\cong A_4$ observe that $V(t_4)\cup\{\infty\}$ and $V(p_4)\cup V(q_4)=V(s_4)$ are already orbits of~$S_4$. 
To obtain the reduced automorphism group $A_4$ the branch locus must therefore contain exactly one of the orbits $V(p_4)$ and $V(q_4)$. As the action of $S_4\setminus A_4$ interchanges these, to avoid duplicity of isomorphic curves we can restrict ourselves to the case that $V(p_4)\subset\Branch(\pi)$. Since $|V(p_4)|=4<6$, this leaves only the case $\Branch(\pi)= V(p_4)\cup V(t_4)\cup\{\infty\}$. 
Then $|\Branch(\pi)|=10$, which does not occur for any larger group; hence the reduced automorphism group of $X$ is really $A_4$ in this case.

For the case $\oG\cong D_n$ observe that $\{0,\infty\}$ and $V(x^n-1)\cup V(x^n+1) = V(x^{2n}-1)$ are already orbits of~$D_{2n}$. To obtain the reduced automorphism group $D_n$ the branch locus must therefore contain exactly one of $V(x^n-1)$ and $V(x^n+1)$. As the action of $D_{2n}\setminus D_n$ interchanges these, to avoid duplicity we can restrict ourselves to the case that $V(x^n-1)\subset\Branch(\pi)$. 
Then $\Branch(\pi)=V(x^n-1)$, respectively $V(x^n-1)\cup\{0,\infty\}$. The condition $|\Branch(\pi)|=2g+2\ge6$ then implies that $n$ is even and at least~$6$, respectively~$4$. If $n=4$ we have $\Branch(\pi)=V(x^4-1)\cup\{0,\infty\}=V(t_4)\cup\{\infty\}$, in which case the reduced automorphism group is~$S_4$, as seen above. Thus we must have $n\ge 6$. Then $D_n$ does not embed into $S_4$ or~$A_5$ and $|\Branch(\pi)|$ is too small to be invariant under a larger dihedral group, so the reduced automorphism group is indeed~$D_n$. This gives precisely two curves for every even $n\ge6$. 

Consider now the case $\oG\cong C_n$. Then $\oG$ acts on $\BP^1_\BC$ by multiplication with $n$-th roots of unity, and $\Branch(\pi)$ must contain precisely one free $\oG$-orbit. This orbit must have the form $V(x^n-a^n)$ for some $a\in\BC^\times$. After rescaling by $x\mapsto ax$, which commutes with the action of $\oG$, we may assume that this orbit is $V(x^n-1)$. If $n$ is even, the parity requirement implies that $\Branch(\pi) = V(x^n-1)$ or $V(x^n-1)\cup\{0,\infty\}$. In both these cases $\Branch(\pi)$ is also invariant under the substitution $x\mapsto x^{-1}$, so that the reduced automorphism group of $X$ contains the dihedral group $D_n$, which is already covered by the preceding case. If $n$ is odd, the parity requirement implies that $\Branch(\pi) = V(x^n-1)\cup\{0\}$ or $V(x^n-1)\cup\{\infty\}$. These cases correspond to each other under the substitution $x\mapsto x^{-1}$ which normalizes~$\oG$, so it suffices to consider the case $V(x^n-1)\cup\{\infty\}$. 
The condition $n+1=|\Branch(\pi)|\ge6$ then requires that $n\ge5$.
We claim that in this situation the reduced automorphism group is really~$C_n$. Indeed, the equality $|\Branch(\pi)|=n+1$ admits no larger cyclic group, and by the preceding case it can admit at most a dihedral group $D_m$ with $|\Branch(\pi)|=m$ or $m+2$. In that case we would have $m=n\pm1$, which renders an embedding $C_n\into D_m$ impossible for $n\ge5$. Of the remaining candidates $A_4$, $S_4$, $A_5$, only one contains a cyclic subgroup of order $n\ge5$, namely $A_5$ with $n=5$. But in each of the cases above with $\oG=A_5$ we have $|\Branch(\pi)|\ge12>5+1$, so $A_5$ cannot occur either. This proves the claim. 

\medskip
All in all this shows that the hyperelliptic curves with many automorphisms are, up to isomorphism, precisely those in Table~\ref{TableHyperellipticManyAuto} of the introduction and that those are pairwise non-isomorphic. 

\begin{Rem}
\rm The curve $X_4$ can already be defined over $\BQ$. Wolfart \cite[\S6.4]{Wolfart2000} gives the affine equation $y^2 = x(3x^4+1)(3x^4+6x^2-1)$ and states that the hyperelliptic curve defined by this equation has many automorphisms and its automorphism group is isomorphic to $\SL_2(3)$. It defines a curve isomorphic to $X_4$ because there is, up to isomorphism, only one hyperelliptic curve of genus 4 with these properties.
\end{Rem}

	
\section{Representation theoretic condition for complex multiplication}
\label{Rep}
	
For the moment let $X$ be an arbitrary smooth complex algebraic curve of genus $g\ge 2$ with automorphism group~$G$. Let $\chi_X$ denote the character of the representation of $G$ on $H^0(X,\Omega_X)$. Let $\Sym^2\chi_X$ denote its symmetric square, that is, the character of $G$ on $\Sym^2_\BC H^0(X,\Omega_X)$, and let $\chi_\triv$ denote the character of the trivial representation on~$\BC$. Then by Streit \cite[Prop.\;5]{Streit2001} we have:

\begin{Prop}\label{RepConditionStreit}
If $\langle \Sym^2\chi_X, \chi_\triv\rangle=0$, then $\Jac(X)$ has complex multiplication.
\end{Prop}

In fact $\Sym^2_\BC H^0(X,\Omega_X)$ is naturally isomorphic to the tangent space of the Siegel moduli space $\CA_g$ at the point corresponding to~$\Jac(X)$, and the assumption implies that $\Jac(X)$ cannot be deformed non-trivially as a polarized abelian variety together with the action of~$G$. The point on $\CA_g$ is therefore an isolated point of a certain Shimura subvariety of PEL-type and hence a special point. By the Torelli theorem $X$ itself can then also not be deformed non-trivially together with the action of~$G$, so $X$ has many automorphisms. 

\medskip
Let us now return to a hyperelliptic curve $X$ as above. To use Streit's condition, we need to calculate $\Sym^2\chi_X$. 

\begin{Prop}\label{RepChiX}
Let $\tau\in G$ and let $n$ be the order of its image $\otau\in\oG$. 
Let $P\in\BP^1_\BC$ be a fixed point of $\otau$ and let $\zeta$ be the eigenvalue of $\otau$ on the tangent space at~$P$. Set $k:=1$ if $P\in\Branch(\pi)$ and $k:=0$ otherwise.
Then we have the following character values:
\[\begin{array}{|c||c|c|c|}
\hline
{\large\strut} \text{Case} 			& n=1 & n=2  & n>2 \\
\hline\hline
{\Large\strut} \chi_X(\tau)			& \pm g & \pm (-1)^{\frac{k}{2}}\frac{(-1)^g-1}{2} & \pm \zeta^{1-\frac{k}{2}}\cdot\frac{\zeta^g-1}{\zeta-1} \\[6pt] \hline
{\Large\strut} \chi_X(\tau^2)		& g & (-1)^k g & \zeta^{2-k}\cdot\frac{\zeta^{2g}-1}{\zeta^2-1} \\[6pt] \hline
{\Large\strut} \Sym^2\chi_X(\tau)	& \frac{g(g+1)}{2} & (-1)^k\cdot \frac{1+(-1)^{g+1}+2g}{4} & \zeta^{2-k}\cdot \frac{(\zeta^g-1)(\zeta^{g+1}-1)}{(\zeta-1)(\zeta^2-1)} \\[6pt] \hline
\end{array}\]
\end{Prop}

\begin{Proof}
After a change of coordinates in $\BP^1_\BC$ mapping $P$ to $0$ and the other fixed point of $\otau$ to $\infty$, we may assume without loss of generality that $\otau$ has the form $x\mapsto \zeta x$. In these coordinates the affine equation for $X$ takes the form $y^2=x^kh(x^n)$ for a separable polynomial $h\in\BC[x]$ with $h(0)\not=0$. The action of $\tau$ is therefore given by
\[\tau: (x,y)\mapsto (\zeta x, \pm \zeta^{\frac{k}{2}}y)\]
for some sign and some choice of the fractional power $\zeta^{\frac{k}{2}}$.
It it well known that the differentials $\frac{x^jdx}{y}$ for all $0\le j\le g-1$ form a basis of $H^0(X,\Omega_X)$. As they are also eigenvectors with eigenvalues $\pm\zeta^{j+1-\frac{k}{2}}$ under~$\tau$, a quick computation now yields the desired values of $\chi_X(\tau)=\trace\tau|_{H^0(X,\Omega_X)}$ and $\chi_X(\tau^2)$.
The values of $\Sym^2\chi_X$ follow from these using the formula $\Sym^2\chi_X(\tau)=\frac{1}{2}(\chi_X(\tau)^2+\chi_X(\tau^2))$ and the fact that $\zeta=-1$ if $n=2$.
\end{Proof}

\begin{Rem}
\rm Since the hyperelliptic involution $\sigma$ acts as multiplication by $-1$ on $H^0(X,\Omega_X)$, it acts trivially on the symmetric square; hence $\Sym^2\chi_X$ comes from a character of~$\oG$. 
The formula for $\Sym^2\chi_X(\tau)$ in Proposition \ref{RepChiX} also depends only on~$\bar\tau$. To compute $\Sym^2\chi_X$ it therefore suffices to work only with~$\oG$.
\end{Rem}

\begin{Prop}
The jacobians of $X_4$, $X_5$, $X_7$, $X_9$ and $X_{14}$ have complex multiplication.
\end{Prop}

\begin{Proof}
By Proposition \ref{RepConditionStreit} it suffices to show that $\langle \Sym^2\chi_X, \chi_\triv\rangle=0$ for each of these curves. We verify this using the formulas in Proposition \ref{RepChiX} and the calculations in the respective worksheet \cite{MuellerPinkWorksheets} using GAP. For $X_4$ and $X_5$ this is already known: Wolfart \cite[\S6]{Wolfart2000} states that the curves $X_4$ and $X_5$ have jacobians with complex multiplication.
\end{Proof}


\section{Quotients of small genus}
\label{Quots}

To show that $\Jac(X)$ does not have complex multiplication, it suffices to exhibit an abelian subvariety without complex multiplication. A natural candidate for this is the jacobian of the quotient of $X$ by a subgroup $H<G$ whose genus is positive but small. We therefore make a list of suitable quotients for the 10 curves left over from the preceding section.

Note that $H$ must not contain the hyperelliptic involution~$\sigma$, because in that case the quotient morphism $X\onto H\backslash X$ factors through $\pi\colon X\onto\BP^1_\BC$, forcing $H\backslash X$ to have genus~$0$. Thus $H$ must map isomorphically to a subgroup $\oH<\oG$. While this $\oH$ does not always determine $H$ up to conjugation, in all our cases we can actually find an equation for $H\backslash X$ by working only with~$\oH$. 

\begin{Prop}\label{QuotsList}
In each of the following cases there is a subgroup $H<G$ mapping isomorphically to $\oH<\oG$ where the quotient $H\backslash X$ of genus $\bar g$ is described in the following table. 
\[\begin{array}{|c||c|c|c|c|}
\hline
{\large\strut} X      & \oG		& \oH		& \bar g		&\text{Affine equation of $H\backslash X$} \\
\hline\hline
{\large\strut} X_6    & S_4		& C_2		& 1				& v^2=u^4 + 14u^2 + 1 		\\ \hline
{\large\strut} X_8    & S_4		& S_3		& 1				& 
v^2=u(u-1)(4u+1+i\sqrt{2})
\\ \hline
{\large\strut} X_{10} & S_4		& C_4		& 2				& v^2=u^5 - 19u^4 - 494u^3 - 494u^2 - 19u + 1 			\\ \hline
{\large\strut} X_{11} & S_4		& C_3		& 4				& 
v^2= u(u-1)(u+8)(u^2+8)(u^2+4u-8)(u^2+8u-8)
\\ \hline
{\large\strut} X_{12} & A_5		& C_5		& 1				& v^2=u(u^2+11u-1)			\\ \hline
{\large\strut} X_{13} & A_5		& C_5		& 1				& v^2=\bar r_5(u) 			\\ \hline
{\large\strut} X_{15} & A_5		& A_4		& 1				& 
v^2=u^3+5u^2+40u
			\\ \hline
{\large\strut} X_{16} & A_5		& C_5		& 4				& v^2=u(u^2 + 11u - 1)\bar t_5(u) 			\\ \hline
{\large\strut} X_{17} & A_5		& C_5		& 4				& v^2=\bar r_5(u)\bar t_5(u)			\\ \hline
{\large\strut} X_{18} & A_5		& C_5		& 6				& v^2=u(u^2 + 11u - 1)\bar r_5(u)\bar t_5(u)			\\ \hline
\end{array}\]
Here $\bar r_5$ and $\bar t_5$ are the unique polynomials satisfying $\bar r_5(x^5)=r_5(x)$ and $\bar t_5(x^5)=\nolinebreak t_5(x)$, so that 
$$\begin{array}{rl}
\bar r_5 &=\ u^4 - 228  u^3 + 494  u^2 + 228  u + 1, \\[3pt]
\bar t_5 &=\ u^6 + 522  u^5 - 10005 u^4 - 10005 u^2 - 522  u + 1.
\end{array}$$
\end{Prop}

\begin{Proof}
In several cases it suffices to note that the given affine equation for $X$ has the form $y^2=g(x^n)$ for some polynomial $g$ and some integer $n>1$, so that $y^2=g(z)$ is the affine equation of a suitable quotient. 
This and more was already done, e.g. for $X_6$ by Wolfart \cite[\S6]{Wolfart2000} and for $X_{12}$ by Paulhus \cite[Thm.\;2]{Paulhus2013}.)

In the general case we proceeded as follows.
For each curve $X$ we know the precise automorphism group $G$ by Table \ref{TableHyperellipticManyAuto}. After a preliminary search through subgroups of $G$ we identified a likely candidate for~$H$. The details of these calculations are immaterial for the desired result, because  we only need to know or guess a good choice of~$\oH$.
So suppose that $\oH$ is given. Let $y^2=f(x)$ be an affine equation for~$X$, where $f$ is a separable polynomial of degree~$n$. 

Assume first that $\oH$ is cyclic of order $n>1$. Then for a suitable fractional linear transformation $x=\frac{az+b}{cz+d}$ the group $\oH$ acts on $z$ by multiplication with $n$-th roots of unity. After the substitution $y=v/(cz+d)^{\lceil\frac{n}{2}\rceil}$ we compute a new affine equation for $X$ in the form $v^2 = z^k g(z^n)$ for a separable polynomial $g$ and an integer $k\in\{0,1\}$. 
If $n$ is odd and $k=0$, the quotient $H\backslash X$ must have the affine equation $v^2 = g(u)$. If $n$ is odd and $k=1$, the further substitution $v=w/z^{\lfloor\frac{n}{2}\rfloor}$ brings the equation into the form $w^2 = z^n g(z^n)$, so the quotient has the affine equation $w^2 = ug(u)$.
If $n$ is even and $k=0$, the formula $v^2 = g(z^n)$ does not tell us whether $v$ is fixed by $H$ or not. If not, the new variable $w:=vz^{\frac{n}{2}}$ is fixed by $H$ and satisfies the equation 
$w^2 = z^ng(z^n)$. The affine equation for $H\backslash X$ is thus either $v^2 = g(u)$ or $w^2 = ug(u)$. The other equation then describes the quotient $H'\backslash X$ for another subgroup $H'<G$ lifting~$\oH$, but if we do not care about the precise form of~$H$, we can continue with either choice. Finally, if $n$ is even and $k=1$, the equation $v^2 = z g(z^n)$ implies that a generator of $H$ which multiplies $z$ with a primitive $n$-th root of unity must multiply $v$ with a primitive $2n$-th root of unity. The $n$-th power of that generator is then the hyperelliptic involution and lies in~$H$, which we had ruled out before (or this $\oH$ was a bad choice).

If $\oH$ is not cyclic, it is nevertheless solvable in each of our cases. Thus by repeating the same procedure over a composition series of~$\oH$ we can find an explicit equation for the desired quotient in that case, too.
(If at an earlier stage with $n$ even and $k=0$ we choose the wrong quotient for a subgroup of~$\oH$, we detect that later when the equation is no longer invariant under~$\oH$.)
Finally, the genus of $H\backslash X$ can be read off directly from the degree of the affine equation for it.

The actual calculations were carried out with Sage in the respective worksheet \cite{MuellerPinkWorksheets}. In some cases we applied an additional fractional linear transformation at the end to simplify the equation.
\end{Proof}

\begin{Cor}\label{EllQuot}
The jacobians of $X_6$, $X_8$, $X_{12}$, $X_{13}$, and $X_{15}$ do not have complex multiplication.
\end{Cor}

\begin{Proof}
By Proposition \ref{QuotsList} each of these curves possesses a quotient of genus~$1$. By the same worksheet \cite{MuellerPinkWorksheets} as above the respective $j$-invariant comes out to be:
\[\begin{array}{|c||c|}
\hline
{\large\strut} X      &\text{$j$-invariant of $H\backslash X$} \\
\hline\hline
{\large\strut} X_6    & 2^4\cdot 13^3\cdot 3^{-2} 		\\ \hline
{\large\strut} X_8    & 
2^4\cdot (1+i\sqrt{2}\,)^4\cdot (19+6i\sqrt{2}\,)^3\cdot 3^{-6}
\\ \hline
{\large\strut} X_{12} & 2^{14}\cdot 31^3\cdot 5^{-3}	\\ \hline
{\large\strut} X_{13} & 2^{17}\cdot 3^{-2} 				\\ \hline
{\large\strut} X_{15} & 2^2\cdot 19^3\cdot 3^{-3}		\\ \hline
\end{array}\]
In each case the $j$-invariant is not an algebraic integer. Thus the elliptic curve $H\backslash X$ does not have complex multiplication, and hence $\Jac(X)$ does not, either. 
(For $X_6$ this is known: Wolfart \cite[\S6.3]{Wolfart2000} shows that $\Jac(X_6)$ does not have complex multiplication by exhibiting, up to a change of coordinates, the same elliptic quotient.)
\end{Proof}

\medskip
For the remaining 5 curves we have to deal with quotients of genus $2$, $4$, and~$6$, for which no suitable analogue of the $j$-invariant criterion is available. We will return to these curves after discussing a different criterion based on characteristic polynomials of Frobenius.


\section{Computational criterion for complex multiplication}
\label{CMCrit}

Given an abelian variety $A$ over a number field~$K$, one can effectively determine the endomorphism ring $\End_{\oK}(A)$ over an algebraic closure $\oK$ of~$K$, at least in principle. All methods for this use the Tate conjecture for endomorphisms and either specific height estimates based on Faltings \cite{FaltingsFiniteness1984} or from Masser and W\"ustholz \cite{MasserWuestholz1995}, or a day and night algorithm as in Lombardo \cite[\S5]{Lombardo2016} (compare also the proof of Lemma 2 in Achter \cite{Achter2005}). In particular one can effectively decide whether $A$ has complex multiplication over~$\oK$. But these methods require complicated code whose computational feasibility is unclear. 

To keep computations to a minimum, we formulate and implement a simple criterion based solely on characteristic polynomials of Frobenius. Let $\Gamma_\ell<\GL_{2\dim(A)}(\BQ_\ell)$ be the image of the $\ell$-adic Galois representation associated to~$A$. Then the Tate conjecture for endomorphisms implies that $A$ has complex multiplication over $\oK$ if and only if $\Gamma_\ell$ possesses a commutative subgroup of finite index. Our criterion exploits the fact that otherwise $\Gamma_\ell$ contains an open subgroup of a non-trivial semisimple algebraic group over~$\BQ_\ell$, which possesses non-isomorphic maximal tori, and that therefore the characteristic polynomials of Frobenius cannot all split over the same number field. The article with Larsen \cite{LarsenPink1992} was based on essentially the same principle.

\medskip
As a preparation consider any monic polynomial $f\in\BQ[T]$ without multiple factors. Then $f(T)={\prod_{i=1}^r(T-\alpha_i)}$ for pairwise distinct $\alpha_i\in\bar\BQ$, and the factor ring $E_f := \BQ[T]/(f)$ is a finite direct product of number fields. Let $t\in E_f$ denote the residue class of~$T$, so that $f$ is precisely the minimal polynomial of $t$ over~$\BQ$. 
Let $E'_f$ denote the intersection of the subrings $\BQ[t^n]\subset E_f$ for all integers $n\ge1$. Since $E_f$ has finite dimension over~$\BQ$ and $\BQ[t^m]\subset\BQ[t^n]$ whenever $n|m$, we have $E'_f=\BQ[t^{n_0}]$ for some integer $n_0\ge1$. Then $E'_f\cong \BQ[T]/(g)$, where $g\in\BQ[T]$ is the minimal polynomial of $t^{n_0}$ over~$\BQ$.

\begin{Prop}\label{E'f=Ef}
One can effectively find such $n_0$ and $g$ and decide whether $E'_f=E_f$. Moreover we have $E'_f=E_f$ if and only if for all $i\not=j$ with $\alpha_j\not=0$ the ratio $\alpha_i/\alpha_j$ is not a root of unity.
\end{Prop}

\begin{Proof}
By construction $t$ is a semisimple element of~$E_f$. Thus $t^n$ is semisimple and its minimal polynomial over $\BQ$ is the product of the linear factors $T-\alpha_i^n$ with all repetitions removed. In particular $\dim_\BQ\BQ[t^n]$ is the number of distinct values among $\alpha_1^n,\ldots,\alpha_r^n$.

Let $S$ be the set of ratios $\alpha_i/\alpha_j$ for all $i\not=j$ which are well-defined and roots of unity. Let $n_0$ be the least common multiple of the orders of all elements of~$S$. Then for any multiple $n$ of~$n_0$ and any $i\not=j$, we have $\alpha_i^n=\alpha_j^n$ if and only if 
$\alpha_i^{n_0}=\alpha_j^{n_0}$, and hence $\BQ[t^n]=\BQ[t^{n_0}]$. Varying $n$ this implies that $E'_f=\BQ[t^{n_0}]$. In the same way we see that $E'_f=\BQ[t^n]$ if and only if $n_0|n$. Thus $n_0$ is the unique smallest integer $n\ge1$ with $E'_f=\BQ[t^n]$. In particular $E'_f=E_f$ if and only if $n_0=1$.

To determine $n_0$ effectively we first replace $f(T)$ by $f(T)/T$ if possible. Next we compute the polynomial $h(T) := \prod_{i,j=1}^r(T-\alpha_i/\alpha_j) \in \BQ[T]$, say, using symmetric polynomials. Then $n_0$ is the least common multiple of the set of integers $n>1$ such that the $n$-th cyclotomic polynomial $\Phi_n$ divides~$h$. Since there there exist only finitely many cyclotomic polynomials of degree $\phi(n)=\deg\Phi_n\le \deg h$ and we can list them all explicitly, we can effectively determine~$n_0$.

Again using symmetric polynomials we can effectively compute the polynomial $k(T) := {\prod_{i=1}^r(T-\alpha_i^{n_0})} \in \BQ[T]$. Then $k/\gcd(k,\frac{dk}{dT})\in\BQ[T]$ is a polynomial with the same roots but all repeating factors removed; hence it is the desired polynomial~$g$.
\end{Proof}

\medskip
Now we fix an abelian variety $A$ of dimension $d$ over a number field~$K$. We call a prime $\Fp$ of the ring of integers of $K$ \emph{good} if $A$ has good reduction at~$\Fp$. For any good prime $\Fp$ let $f_\Fp\in\BQ[T]$ denote the minimal polynomial of the Frobenius endomorphism $\Frob_\Fp$ acting on the Tate module of the reduction of~$A$. As this endomorphism is semisimple, the polynomial $f_\Fp$ has no multiple factors, and we can apply the preceding theory to it. We call a good prime $\Fp$ \emph{very good} if $E'_{f_\Fp}=E_{f_\Fp}$ in the notation above.

\begin{Thm}\label{CMCritThm}
The following are equivalent:
\begin{enumerate}
\item[(a)] The abelian variety $A$ has complex multiplication over~$\oK$.
\item[(b)] There exists a commutative semisimple $\BQ$-algebra $E$ of dimension $\le 2d$ such that for every good prime $\Fp$ there exists an embedding $E'_{f_\Fp}\into E$.
\item[(c)] There exists a commutative semisimple $\BQ$-algebra $E$ of dimension $\le 2d$ such that for every very good prime $\Fp$ there exists an embedding $E_{f_\Fp}\into E$.
\item[(d)] There exists a number field $F$ such that for every very good prime $\Fp$ the polynomial $f_\Fp$ splits completely over~$F$.
\item[(e)] There exists a prime number $\ell$ such that for every very good prime $\Fp$ the polynomial $f_\Fp$ splits completely over~$\BQ_\ell$.
\end{enumerate}
\end{Thm}

\begin{Proof}
We first show that (a) implies (b). So assume that $A$ has complex multiplication over~$\oK$. By definition this means that $R := \End_{\oK}^\circ(A)$ contains a commutative semisimple $\BQ$-algebra $E$ of dimension $2d$. We claim that $E$ has the property in (b). 
To see this consider any good prime~$\Fp$. 
Note that $R=\End_{K'}^\circ(A)$ for some finite extension $K'$ of~$K$.
Let $\Fp'$ be a prime of the ring of integers of $K'$ that lies above~$\Fp$, and let $k_\Fp\subset k_{\Fp'}$ denote the corresponding residue fields. Then the reduction of endomorphisms induces an injective homomorphism $R\into R_{\Fp'} := \End^\circ_{k_{\Fp'}}(A_\Fp)$. The image of $E$ in $R_{\Fp'}$ is thus a commutative semisimple $\BQ$-subalgebra of dimension $2d$. By Chai, Conrad and Oort \cite[Thm.\;1.3.1.1]{ChaiConradOort2014}, this image is therefore its own commutant in~$R_{\Fp'}$. It follows that the image of $E$ contains the center of $R_{\Fp'}$. As the Frobenius endomorphism $\Frob_{\Fp'}$ lies in the center of $R_{\Fp'}$, this yields an embedding of $\BQ$-algebras $\BQ[\Frob_{\Fp'}]\into E$.

Now observe that the minimal polynomial $f_\Fp$ of $\Frob_\Fp$ as an endomorphism of $V_\ell(A_\Fp)$ is also its minimal polynomial over $\BQ$ as an element of the $\BQ$-algebra~$R_{\Fp'}$.
Since $\Frob_\Fp$ is semisimple, the subalgebra $\BQ[\Frob_\Fp]\subset R_{\Fp'}$ is therefore isomorphic to $E_{f_\Fp}:=\BQ[T]/(f_\Fp)$ such that $\Frob_\Fp$ corresponds to the residue class 
$t\in E_{f_\Fp}$ of~$T$.
On the other hand we have $\Frob_{\Fp'}=\Frob_\Fp^n$, where $n$ denotes the degree of the residue field extension $k_{\Fp'}/k_\Fp$. The subalgebra $\BQ[\Frob_{\Fp'}]\subset R_{\Fp'}$ is therefore isomorphic to the subalgebra $\BQ[t^n]\subset E_{f_\Fp}$. By the definition of $E'_{f_\Fp}$ there thus exists an embedding of $\BQ$-algebras $E'_{f_\Fp} \into \BQ[t^n] \cong \BQ[\Frob_{\Fp'}]$. Combined with the embedding $\BQ[\Frob_{\Fp'}]\into E$, this yields the desired embedding $E'_{f_\Fp} \into E$. This finishes the proof that (a) implies (b).

\medskip
Clearly (b) implies (c). Next, given any $E$ as in (c), choose a number field $F$ that is Galois over $\BQ$ such that every simple factor of $E$ embeds into~$F$. Then by the construction of $E_{f_\Fp}$, for every very good prime $\Fp$ the polynomial $f_\Fp$ splits completely over~$F$. Thus (c) implies (d). Also (d) implies (e) by taking any prime number $\ell$ that splits completely in~$F$.

\medskip
It remains to show that (e) implies (a). So let $\ell$ be a prime number as in (e). Consider the associated Galois representation $\rho_\ell\colon \Gal(\oK/K)\to\Aut_{\BQ_\ell}(V_\ell(A))\cong\GL_{2d}(\BQ_\ell)$, let $\Gamma_\ell$ denote its image, and let $G_\ell$ be the Zariski closure of $\Gamma_\ell$ in the algebraic group $\GL_{2d,\BQ_\ell}$. By Faltings \cite{FaltingsFiniteness1984} the representation $\rho_\ell$ is semisimple and $\End_\oK^\circ(A)\otimes_\BQ\BQ_\ell \cong \End_{\BQ_\ell[G_\ell^\circ]}(V_\ell(A))$. Thus the identity component $G_\ell^\circ$ is a reductive linear algebraic group. If $G_\ell^\circ$ is a torus, it follows that $\End_\oK^\circ(A)\otimes_\BQ\BQ_\ell$ contains a commutative semisimple $\BQ_\ell$-subalgebra of dimension~$2d$. Since $\End_\oK^\circ(A)$ is a semisimple $\BQ$-algebra, it then contains a commutative semisimple $\BQ$-sub\-alge\-bra of dimension~$2d$; so $A$ has complex multiplication over~$\oK$. 

Suppose that $A$ does not have complex multiplication over~$\oK$. Then $G_\ell^\circ$ is not a torus; hence its derived group $(G_\ell^\circ)^\der$ is non-trivial. Then $(G_\ell^\circ)^\der$ possesses a non-split maximal torus defined over~$\BQ_\ell$, for instance because every maximal torus is non-split if $(G_\ell^\circ)^\der$ is non-split and by \cite[Thm.\;3.6]{LarsenPink1992} otherwise. Thus $G_\ell^\circ$ possesses a non-split maximal torus~$T_\ell$.

Choose a finite extension $K'\subset\oK$ of $K$ such that $\Gamma_\ell' := \rho_\ell(\Gal(\oK/K'))$ is contained in $G_\ell^\circ(\BQ_\ell)$ and consists of matrices that are congruent to the identity modulo~$\ell$. For later use we fix an auxiliary prime $\ell'\not=\ell$ and, after enlarging~$K'$, assume that the $\ell'$-adic representation of $\Gal(\oK/K')$ is also congruent to the identity modulo~$\ell'$.

By Bogomolov \cite{Bogomolov1980} the subgroup $\Gamma'_\ell<G_\ell(\BQ_\ell)$ is open; hence it contains an open subgroup of $T_\ell(\BQ_\ell)$; so there exists a regular semisimple element $t\in\Gamma'_\ell\cap T_\ell(\BQ_\ell)$. Thus $T_\ell$ is the centralizer of $t$ in $G_\ell^\circ$. As $T_\ell$ is non-split, it follows that the characteristic polynomial of $t$ does not split over~$\BQ_\ell$. 
Now observe that the set of monic polynomials of degree $2d$ over $\BQ_\ell$ which are split over $\BQ_\ell$ is the image of a proper continuous map $\BQ_\ell^d\to\BQ_\ell[T]$, $(\alpha_1,\ldots,\alpha_{2d})\mapsto \prod_{i=1}^{2d}(T-\alpha_i)$. It is thus a closed subset of the set of all monic polynomials of degree $2d$ over~$\BQ_\ell$. 
Since it does not contain the characteristic polynomial of~$t$, and the characteristic polynomial of $g\in\GL_{2d}(\BQ_\ell)$ 
varies continuously with~$g$, 
there exists an open neighborhood $U\subset\GL_{2d}(\BQ_\ell)$ of $t$ such that the characteristic polynomial of any $g\in U$ does not split over~$\BQ_\ell$. Note that $\Gamma'_\ell\cap U$ is then a non-empty open subset of~$\Gamma'_\ell$.

Now recall that the primes $\Fp'$ of absolute degree $1$ of the ring of integers of~$K'$ form a subset of Dirichlet density~$1$. This remains true if we restrict ourselves to primes that do not divide~$\ell$ or any prime of $K$ where $A$ has bad reduction. By the Cebotarev density theorem there therefore exists a prime $\Fp'$ with these properties, such that $\gamma := \rho_\ell(\Frob_{\Fp'})$ lies in $\Gamma'_\ell\cap U$. By construction the prime $\Fp$ of $K$ below $\Fp'$ is then good. Also, the fact that $\Fp'$ has absolute degree $1$ implies that the residue field extension $k_{\Fp'}/k_\Fp$ is trivial; hence $\gamma = \rho_\ell(\Frob_\Fp)$ as well. 
Then the minimal polynomial of $\gamma$ in the given representation is precisely the polynomial $f_\Fp$ from before. Write $f_\Fp(T)=\prod_{i=1}^r(T-\alpha_i)$ with pairwise distinct $\alpha_i\in\bar\BQ$. Then by the choice of $K'$ we have $f_\Fp(T)\equiv (T-1)^r$ modulo $\ell$ and modulo~$\ell'$. Thus for any $i\not=j$ the algebraic number $\alpha_i/\alpha_j$ is congruent to $1$ modulo a prime above $\ell$ and modulo a prime above $\ell'$. Since $\ell$ and $\ell'$ are distinct, this implies that $\alpha_i/\alpha_j$ cannot be a non-trivial root of unity. By Proposition \ref{E'f=Ef} it follows that $\Fp$ is a very good prime. 

Finally, since $\gamma$ lies in~$U$, the characteristic polynomial and hence also the minimal polynomial of $\gamma$ does not split over~$\BQ_\ell$. As this minimal polynomial is~$f_\Fp$, we have achieved a contradiction to (e). 
Therefore (e) implies (a), as desired.
\end{Proof}


\begin{Rem}\label{CMCritRem3}
\rm If $A$ has no complex multiplication over~$\oK$, Theorem \ref{CMCritThm} provides a good chance to prove that efficiently. But it cannot prove that $A$ has complex multiplication with a finite calculation. Yet, in that case the computation gives some information about the endomorphism ring which might indicate a direction where one can find endomorphisms.
\end{Rem}

\begin{Rem}\label{CMCritRem1}
\rm In the proof of the implication (c)$\Rightarrow$(d) of Theorem \ref{CMCritThm}, it follows from a result of Dodson \cite[Prop.\;1.1]{Dodson1984} about the structure of Galois groups of CM-fields that the number field $F$ has degree at most $2^dd!$ over~$\BQ$. Thus if one has found very good primes $\Fp_1,\dots,\Fp_r$ such that the splitting field of the product $f_{\Fp_1}\cdots f_{\Fp_r}$ has degree $>2^dd!$ over~$\BQ$, one can deduce that $A$ has no complex multiplication over~$\oK$. But while this provides a simple algorithm in principle, the size of $2^dd!$ made it impracticable in some of our cases.
\end{Rem}

\begin{Rem}\label{CMCritRem2}
\rm If one knows for some reason that $A$ is simple over~$\oK$, the ring $R := \End_{K'}^\circ(A)$ in the proof of Theorem \ref{CMCritThm} is simple; hence its center $E$ is a field. In that case we obtain the same equivalences with a field $E$ in (b) and (c). But this observation presumes information about the endomorphism ring which we might not have in a given situation. 
Yet perhaps one finds a very good prime $\Fp$ such that the characteristic polynomial of $\Frob_\Fp$ is already irreducible over~$\BQ$. Then the reduction $A_\Fp$ is simple over $\bar k_\Fp$; hence $A$ is simple over~$\oK$.
\end{Rem}

If the rings $E'_{f_\Fp}$ for good primes $\Fp$ are not fields, it may be messy to juggle with their simple factors. At least we can say:

\begin{Cor}\label{CMCritCor1}
If $A$ has complex multiplication over~$\oK$, for any good primes $\Fp_1,\dots,\Fp_r$ the tensor product $E'_{f_{\Fp_1}}\!\otimes_\BQ\ldots\otimes_\BQ E'_{f_{\Fp_r}}$ possesses a simple direct summand of dimension $\le 2d$.
\end{Cor}

\begin{Proof}
The embeddings $E'_{f_{\Fp_i}}\!\into E$ from Theorem \ref{CMCritThm} (b) induce a unitary and hence non-zero ring homomorphism $E'_{f_{\Fp_1}}\! \otimes_\BQ\ldots\otimes_\BQ E'_{f_{\Fp_r}}\! \to E$, and any simple direct summand of its image has dimension $\le 2d$.
\end{Proof}

\medskip
In the present paper we only need the following special case:

\begin{Cor}\label{CMCritCor2}
If $A$ has complex multiplication over~$\oK$, for any very good primes $\Fp_1,\ldots,\Fp_r$ such that the minimal polynomials $f_{\Fp_i}$ are irreducible over~$\BQ$ and the associated fields $E_{f_{\Fp_i}}$ are linearly disjoint over~$\BQ$, we have
$$\prod_{i=1}^r [E_{f_{\Fp_i}}:\BQ]\ \le\ 2d.$$
\end{Cor}

\begin{Proof}
In this case each $E'_{f_{\Fp_i}} = E_{f_{\Fp_i}}$ is a field (which means that $A_{\Fp_i}$ is isotypic), and so is their tensor product over~$\BQ$; hence the inequality follows directly from Corollary \ref{CMCritCor1}.
\end{Proof}


\section{Application of the criterion}
\label{Appl}

\begin{Prop}\label{HyperEllQuot}
The jacobians of $X_{10}$, $X_{11}$, $X_{16}$, $X_{17}$, and $X_{18}$ do not have complex multiplication.
\end{Prop}

\begin{Proof}
For each of these curves $X$ the quotient $\oX := H\backslash X$ from Proposition \ref{QuotsList} can be defined over~$\BQ$. Let $\bar g$ be its genus and $y^2=f(x)$ an affine equation for it with a separable polynomial $f\in\BZ[x]$. For any prime number $p>2$ that does not divide the discriminant of~$f$, the curve $\oX$ has good reduction at~$p$. Then its jacobian $\oA$ also has good reduction at~$p$, and the reduction $\oA_p$ of $\oA$ is naturally isomorphic to the jacobian of the reduction $\oX_p$ of~$\oX$, for instance by \cite[Prop.\;9.5.20]{FGAexplained}.
The computer algebra system Sage provides a function, which partially uses PARI \cite{PARI2}, that efficiently computes the characteristic polynomial of Frobenius for any hyperelliptic curve over a finite field of odd characteristic. 
Thus for any good prime $p$ as above, we can determine the characteristic polynomial $g_p$ of Frobenius for~$\oA_p$. 
By a direct calculation using Proposition \ref{E'f=Ef} we check if $p$ is very good, that is, if none of the ratios of the roots of $g_p$ is a nontrivial root of unity. By factoring $g_p$ we then obtain the minimal polynomial $f_p$ of $\Frob_p$ and know whether that is irreducible over~$\BQ$.

Repeating this process for sufficiently many good primes, in each case we actually find very good primes $p_1,\dots,p_r$ for which the $f_{p_i}$ are irreducible over~$\BQ$, the associated fields $E_{f_{p_i}}$ are linearly disjoint over~$\BQ$, and 
\[\prod_{i=1}^r\; [E_{f_{p_i}}:\BQ]\ >\ 2\bar g.\]
By Corollary \ref{CMCritCor2} it follows that $\oA$ does not have complex multiplication, and hence $\Jac(X)$ does not, either.
The following table lists the very good primes that we employed:
\[\begin{array}{|c||c|c|}
\hline
{\large\strut} X      & \bar g	&\text{very good primes}\\ \hline\hline
{\large\strut} X_{10} & 2		& 37, 61, 157\\ \hline
{\large\strut} X_{11} & 4		& 7, 73\\ \hline
{\large\strut} X_{16} & 4		& 31, 151\\ \hline
{\large\strut} X_{17} & 4		& 31, 41\\ \hline
{\large\strut} X_{18} & 6		& 131, 211\\ \hline
\end{array}\]
For the actual calculations see the respective worksheet in  \cite{MuellerPinkWorksheets}.
\end{Proof}



\end{document}